\documentclass[a4paper,12pt,twoside]{article}
\usepackage{amsmath,amssymb,ifthen}
\usepackage[amsmath,thmmarks]{ntheorem}
\usepackage{paper}
\usepackage[all]{xy}
\usepackage{mathrsfs}
\SelectTips{cm}{11}

\DeclareMathOperator{\Ind}{Ind}
\DeclareMathOperator{\JL}{\mathit{JL}}
\DeclareMathOperator{\LJ}{\mathit{LJ}}
\DeclareMathOperator{\Irr}{\mathbf{Irr}}
\DeclareMathOperator{\Disc}{\mathbf{Disc}}
\DeclareMathOperator{\Cusp}{\mathbf{Cusp}}
\DeclareMathOperator{\Nilp}{\mathbf{Nilp}}
\DeclareMathOperator{\Nrd}{Nrd}
\DeclareMathOperator{\Rk}{Rk}
\DeclareMathOperator{\Fix}{Fix}
\DeclareMathOperator{\vol}{vol}
\DeclareMathOperator{\cInd}{c-Ind}

\newcommand{\Dr}{\mathrm{Dr}}

\allowdisplaybreaks

\title{Geometric approach to the local Jacquet-Langlands correspondence}
\author{Yoichi Mieda}
\def\headertitle{Geometric approach to the local Jacquet-Langlands correspondence}

\begin{document}

\maketitle

\begin{firstfootnote}
 Faculty of Mathematics, Kyushu University, 744 Motooka, Nishi-ku, Fukuoka, 819--0395 Japan

 E-mail address: \texttt{mieda@math.kyushu-u.ac.jp}

\end{firstfootnote}

\begin{abstract}
 In this paper, we give a purely geometric approach to the local Jacquet-Langlands correspondence
 for $\GL(n)$ over a $p$-adic field, under the assumption that 
 the invariant of the division algebra is $1/n$.
 We use the $\ell$-adic \'etale cohomology of the Drinfeld tower
 to construct the correspondence at the level of the Grothendieck groups with rational coefficients.
 Moreover, assuming that $n$ is prime, we prove that this correspondence
 preserves irreducible representations.
 This gives a purely local proof of the local Jacquet-Langlands correspondence
 in this case.
 We need neither a global automorphic technique nor detailed
 classification of supercuspidal representations of $\GL(n)$. 
\end{abstract}

\section{Introduction}
Let $F$ be a $p$-adic field, i.e., a finite extension of $\Q_p$.
Let $n\ge 1$ be an integer and $D$ a central division algebra over $F$ such that $\dim_FD=n^2$.
The famous local Jacquet-Langlands correspondence gives a natural bijective correspondence between
irreducible discrete series representations of $\GL_n(F)$ and irreducible smooth representations of $D^\times$.
Let us recall its precise statement.
Write $\Irr(D^\times)$ for the set of isomorphism classes of irreducible smooth representations of $D^\times$.
We denote by $\Disc(\GL_n(F))$ the set of isomorphism classes of irreducible discrete series 
representations of $\GL_n(F)$.
For $\rho\in \Irr(D^\times)$ (resp.\ $\pi\in \Disc(\GL_n(F))$), we denote the character of $\rho$ (resp.\ $\pi$)
by $\theta_\rho$ (resp.\ $\theta_\pi$). Here $\theta_\rho$ is a locally constant function on $D^\times$,
and $\theta_\pi$ is a locally integrable function on $\GL_n(F)$ which is locally constant on 
$\GL_n(F)^{\mathrm{reg}}$, the set of regular elements of $\GL_n(F)$.
The precise statement of the local Jacquet-Langlands correspondence is the following:

\begin{thm}[the local Jacquet-Langlands correspondence]\label{thm:LJLC}
 There exists a unique bijection
 \[
  \JL\colon \Irr(D^\times)\yrightarrow{\cong} \Disc\bigl(\GL_n(F)\bigr)
 \]
 satisfying the following character relation:
 for every regular element $h$ of $D^\times$, $\theta_\rho(h)=(-1)^{n-1}\theta_{\JL(\rho)}(g_h)$,
 where $g_h$ is an arbitrary element of $\GL_n(F)$ whose minimal polynomial is the same as that of $h$.
\end{thm}

The original proof of this theorem, due to Deligne-Kazhdan-Vigneras \cite{MR771672} and Rogawski \cite{MR700135},
was accomplished by using a global automorphic method.
In some cases, more explicit local studies can be found in
\cite{MR1263525}, \cite{MR1771941}, \cite{MR2130587}, which are based on the theory of types.
However, apart from the case of $\GL(2)$, a purely local proof of Theorem \ref{thm:LJLC}
seems not to be known yet (\cf \cite[a comment after Theorem 2]{MR2275640}).

In this article, under the assumption that the invariant of $D$ is $1/n$,
we will give a geometric approach to construct the bijection $\JL$ above.
Let $R(D^\times)$ be the Grothendieck group of finite-dimensional smooth representations
of $D^\times$, and $\overline{R}(\GL_n(F))$ the Grothendieck group of finite length
smooth representations of $\GL_n(F)$ ``modulo induced representations''
(for a precise definition, see \cite{MR874042}).
It is known that the classes of elements of $\Irr(D^\times)$ (resp.\ $\Disc(\GL_n(F))$)
form a basis of $R(D^\times)$ (resp.\ $\overline{R}(\GL_n(F))$).
Put $R(D^\times)_\Q=R(D^\times)\otimes_\Z\Q$ and 
$\overline{R}(\GL_n(F))_\Q=\overline{R}(\GL_n(F))\otimes_\Z\Q$.
The main theorems of this article are the following:

\begin{thm}[Theorem \ref{thm:main}]\label{thm:main-intro}
 We can construct the following two homomorphisms geometrically:
 \[
  \JL\colon R(D^\times)_\Q\longrightarrow \overline{R}\bigl(\GL_n(F)\bigr)_\Q,
 \quad
 \LJ\colon \overline{R}\bigl(\GL_n(F)\bigr)_\Q\longrightarrow R(D^\times)_\Q.
 \]
 These two maps are inverse to each other, and satisfy
 the character relations
 \[
  \theta_\rho(h)=(-1)^{n-1}\theta_{\JL(\rho)}(g_h),\quad 
 \theta_\pi(g_h)=(-1)^{n-1}\theta_{\LJ(\pi)}(h)
 \] 
 for every regular $h\in D^\times$.
\end{thm}

\begin{thm}[Theorem \ref{thm:LJLC-prime}]\label{thm:JL-intro}
 If $n$ is prime, then $\JL$ induces a bijection 
 \[
  \JL\colon \Irr(D^\times)\yrightarrow{\cong} \Disc\bigl(\GL_n(F)\bigr).
 \]
\end{thm}
Theorem \ref{thm:JL-intro} provides a purely local proof of Theorem \ref{thm:LJLC} in the case above.
In particular, the local Jacquet-Langlands correspondence for $\GL_2(F)$ and $\GL_3(F)$ are fully recovered.

The geometric object we use is the Drinfeld tower for $\GL_n(F)$. 
It is a tower of rigid spaces over (a disjoint union of) the $(n-1)$-dimensional Drinfeld upper half space
$\P^{n-1}_F\setminus \bigcup_{H}H$,
where $H$ runs through hyperplanes of $\P^{n-1}_F$ defined over $F$
(for more detailed explanation, see Section \ref{sec:Drinfeld-tower}).
Thanks to extensive studies by many people (\cf \cite{MR1464867}, \cite{MR1876802}, \cite{MR2511742}, \cite{MR2308851}), it is now well-known that the local Jacquet-Langlands correspondence is realized
in the $\ell$-adic cohomology $H_\Dr$ of the Drinfeld tower.
The methods in the works cited above are again global and automorphic.
However, there is also a purely local study of the cohomology due to Faltings \cite{MR1302321}.
He began with an irreducible smooth representation $\rho$ of $D^\times$, and
investigated the $\rho$-isotypic part $H_\Dr[\rho]$ of $H_\Dr$ by means of
the Lefschetz trace formula. By this method, he succeeded to observe that the character relation
in Theorem \ref{thm:LJLC} appears naturally in $H_\Dr[\rho]$. 
By using this result, we can give the map $\JL$ in Theorem \ref{thm:main-intro}.

To construct the inverse map $\LJ$, we need to consider the opposite direction;
we begin with an irreducible discrete series representation $\pi$ of $\GL_n(F)$
and investigate the ``$\pi$-isotypic part'' of $H_\Dr$.
More precisely, we should consider the alternating sum of the extension groups 
$\sum_j(-1)^j\Ext^j(H_\Dr,\pi)$, since $\pi$ is neither projective nor injective in general. 
To study it, we apply the method introduced in \cite{LT-LTF}; namely, we use local harmonic analysis,
such as transfer of orbital integrals.
Furthermore, to prove Theorem \ref{thm:JL-intro}, the non-cuspidality result obtained in \cite{non-cusp}
plays a crucial role.

Since our approach is entirely geometric, it is natural to expect that a similar argument may give
interesting consequences for the mod-$\ell$ Jacquet-Langlands correspondence (\cf \cite{Dat-jlmodl}).
The author also expects that our strategy can be extended to other Rapoport-Zink spaces,
especially the Rapoport-Zink space for $\mathrm{GSp}(4)$.
He hopes to deal with these problems in future works.

We sketch the outline of this paper.
In Section \ref{sec:Drinfeld-tower}, we recall the definition of the Drinfeld tower
and results in \cite{MR1302321}. We use the framework of \cite{MR1779801} to deal with
finitely generated representations systematically.
In Section \ref{sec:harmonic-analysis}, we study the alternating sum of the extention groups
$\sum_j(-1)^j\Ext^j(H_\Dr,\pi)$ by means of local harmonic analysis.
In Section \ref{sec:Faltings-isom}, we apply another deep result of Faltings on the comparison
of the Lubin-Tate tower and the Drinfeld tower. It provides a very important finiteness
result on $H_\Dr$. Under this finiteness, results in Section \ref{sec:Drinfeld-tower} and
Section \ref{sec:harmonic-analysis} can be written in a very simple form.
After short preliminaries on representation theory in Section \ref{sec:rep-theory},
finally in Section \ref{sec:LJLC}, we construct the maps $\JL$ and $\LJ$
in Theorem \ref{thm:main-intro} by using the $\ell$-adic cohomology of the Drinfeld tower,
and investigate their properties.

\bigbreak

\noindent\textbf{Acknowledgment}\quad
The author would like to thank Matthias Strauch for very helpful discussions.
He is also grateful to Tetsushi Ito and Sug Woo Shin for their valuable comments.

\bigbreak

\noindent{\bfseries Notation}\quad
For a totally disconnected locally compact group $H$,
let $\Irr(H)$ be the set of isomorphism classes of irreducible smooth representations of $H$.
We denote the Grothendieck group of finitely generated (resp.\ finite length) smooth $H$-representations
by $K(H)$ (resp.\ $R(H)$). Put $R(H)_\Q=R(H)\otimes_\Z\Q$.
For a finite-dimensional smooth representation $\sigma$ of $H$, write $\theta_\sigma$ for the character
of $\sigma$. It is a locally constant function on $H$.
If moreover a Haar measure on $H$ is fixed,
we denote by $\mathcal{H}(H)$ the Hecke algebra of $H$, namely,
the abelian group of locally constant compactly supported functions on $H$ with convolution product.
Put $\overline{\mathcal{H}}(H)=\mathcal{H}(H)/[\mathcal{H}(H),\mathcal{H}(H)]=\mathcal{H}(H)_H$
(the $H$-coinvariant quotient).

Let $F$ be a $p$-adic field and $\mathcal{O}$ its ring of integers.
We denote the normalized valuation of $F$ by $v_F$ and the cardinality of the residue field of $\mathcal{O}$
by $q$.
Fix a uniformizer $\varpi$ of $\mathcal{O}$.
Denote the completion of the maximal unramified extension of $\mathcal{O}$ by $\breve{\mathcal{O}}$ and
the fraction field of $\breve{\mathcal{O}}$ by $\breve{F}$.

Throughout this paper, we fix an integer $n\ge 1$. 
Let $D$ be the central division algebra over $F$ with invariant $1/n$,
and $\mathcal{O}_D$ its maximal order. Fix a uniformizer $\Pi\in \mathcal{O}_D$ such that $\Pi^n=\varpi$.

For simplicity, put $G=\GL_n(F)$. We denote by $G^{\mathrm{reg}}$ (resp.\ $G^{\mathrm{ell}}$)
the set of regular (resp.\ regular elliptic) elements of $G$. Write $Z_G$ for the center of $G$.
We apply these notations to other groups.
For example, we write $(D^\times)^{\mathrm{reg}}$ for the set of regular elements
of $D^\times$.
As in Theorem \ref{thm:LJLC}, for $h\in (D^\times)^{\mathrm{reg}}$,
let $g_h$ be an element of $G^{\mathrm{ell}}$
whose minimal polynomial is the same as that of $h$.
Such an element always exists, and is unique up to conjugacy.
Moreover, $h\mapsto g_h$ induces a bijection between conjugacy classes in $(D^\times)^{\mathrm{reg}}$ and
those in $G^{\mathrm{ell}}$. Therefore, to $g\in G^{\mathrm{ell}}$ we can attach an element
$h_g\in (D^\times)^{\mathrm{reg}}$ whose minimal polynomial is the same as that of $g$.

For a smooth $G$-representation $\pi$ of finite length, we denote by $\theta_\pi$
the distribution character of $\pi$. It is a locally integrable function on $G$ which is locally constant
on $G^{\mathrm{reg}}$.

We identify $F^\times$ with $Z_G$ and $Z_{D^\times}$. Then, we can consider the quotient groups
$G/\varpi^\Z$ and $D^\times/\varpi^\Z$ under a discrete subgroup $\varpi^\Z$ of $F^\times$.
We regard $\Irr(D^\times/\varpi^\Z)$ (resp.\ $\Irr(G/\varpi^\Z)$) as a subset of
$\Irr(D^\times)$ (resp.\ $\Irr(G)$). 
Similarly, $R(D^\times/\varpi^\Z)$ (resp.\ $R(G/\varpi^\Z)$) is regarded as a submodule of
$R(D^\times)$ (resp.\ $R(G)$).
Write $\mathbf{Disc}(G)$ for the set of isomorphism classes of irreducible discrete series representations
of $G$. Put $\Disc(G/\varpi^\Z)=\Disc(G)\cap \Irr(G/\varpi^\Z)$.

Fix Haar measures on $G$ and $D^\times$. We endow $\varpi^\Z$ with the counting measure and
consider the quotient measures on $G/\varpi^\Z$ and $D^\times/\varpi^\Z$.
For $\varphi\in \mathcal{H}(G/\varpi^\Z)$ and $g\in G^\mathrm{ell}$,
put $O^{G/\varpi^\Z}_g(\varphi)=\int_{G/\varpi^\Z}\varphi(x^{-1}gx)dx$ (the orbital integral).
It is well-known that this integral converges.
Similarly, for $\varphi'\in \mathcal{H}(D^\times/\varpi^\Z)$ and $h\in D^\times$,
put $O^{D^\times/\varpi^\Z}_h(\varphi')=\int_{D^\times/\varpi^\Z}\varphi'(y^{-1}hy)dy$.

For a field $k$, we denote its algebraic closure by $\overline{k}$.
Let $\ell$ be a prime which is invertible in $\mathcal{O}$. We fix an isomorphism $\overline{\Q}_\ell\cong \C$
and identify them. Every representation is considered over $\C$.

\section{Drinfeld tower}\label{sec:Drinfeld-tower}
Let us briefly recall the definition of the Drinfeld tower.
For more detailed description, see \cite{MR0422290}, \cite{MR1141456}, \cite[Chapter 3]{MR1393439}.

First of all, fix a special formal $\mathcal{O}_D$-module $\mathbb{X}$
of $\mathcal{O}_F$-height $n^2$ over $\overline{\F}_q$.
It is well-known that such $\mathbb{X}$ is unique up to $\mathcal{O}_D$-isogeny.

We denote by $\Nilp$ the category of $\breve{\mathcal{O}}$-schemes on which $\varpi$ is locally nilpotent.
Consider the functor $\mathcal{M}_{\Dr}$ from $\Nilp$ to the category of sets that maps $S$ to the set of
isomorphism classes of pairs $(X,\rho)$ consisting of
\begin{itemize}
 \item a special formal $\mathcal{O}_D$-module $X$ over $S$,
 \item and an $\mathcal{O}_D$-quasi-isogeny $\rho\colon \mathbb{X}\otimes_{\overline{\F}_q}\overline{S}\longrightarrow X\otimes_S\overline{S}$,
\end{itemize}
where we put $\overline{S}=S\otimes_{\breve{\mathcal{O}}}\overline{\F}_q$.
Then $\mathcal{M}_{\Dr}$ is represented by a formal scheme locally of finite type over $\breve{\mathcal{O}}$.
We denote the formal scheme by $\mathcal{M}_{\Dr}$ again, and the rigid generic fiber of $\mathcal{M}_{\Dr}$
by $M_{\Dr}$. It is known that $M_{\Dr}$ is the disjoint union of countable copies of
the $(n-1)$-dimensional Drinfeld upper half space.

For an integer $m\ge 0$, let $M_{\Dr,m}$ be the rigid space classifying $\Pi^m$-level structures
on the universal formal $\mathcal{O}_D$-module over $M_\Dr$. It is a finite \'etale Galois covering of $M_\Dr$
with Galois group $(\mathcal{O}_D/(\Pi^m))^\times$. The projective system $\{M_{\Dr,m}\}_{m\ge 0}$ is called
the Drinfeld tower.
We can define a natural right action of $G=\GL_n(F)$ on each $M_{\Dr,m}$ because $G$ is isomorphic to
the group of self $\mathcal{O}_D$-quasi-isogenies of $\mathbb{X}$.
On the other hand, $D^\times$ also acts naturally on $M_{\Dr,m}$ on the right, since $1+\Pi^m\mathcal{O}_D$ is
a normal subgroup of $D^\times$.

Since $M_m$ is too large (it has infinitely many connected components),
we take the quotient $M_{\Dr,m}/\varpi^\Z$ of $M_{\Dr,m}$ by $\varpi^\Z\subset F^\times=Z_G\subset G$.
Put
\[
 H^i_{\Dr,m}=H^i_c\bigl((M_{\Dr,m}/\varpi^\Z)\otimes_{\breve{F}}\overline{\breve{F}},\overline{\Q}_\ell\bigr),\quad 
 H^i_{\Dr}=\varinjlim_m H^i_{\Dr,m}.
\]
These are smooth representations of $G/\varpi^\Z\times D^\times/\varpi^\Z$. Unless $n-1\le i\le 2(n-1)$,
$H^i_{\Dr,m}=H^i_\Dr=0$.

\begin{prop}\label{prop:fin-gen}
 The representation $H^i_{\Dr,m}$ is finitely generated as a $G$-module.
 Moreover, there exist compact open subgroups $K_1,\ldots,K_N$ of $G/\varpi^\Z$, 
 $\varepsilon_\nu\in \{\pm 1\}$ and finite-dimensional smooth representations
 $\sigma_{m,\nu}$ of $K_\nu\times D^\times/\varpi^\Z$
 for each $1\le \nu\le N$
 such that the following holds:
 \[
  \sum_i (-1)^i[H^i_{\Dr,m}]=\sum_{\nu=1}^N \varepsilon_\nu [\cInd_{K_\nu}^{G/\varpi^\Z}\sigma_{m,\nu}]\quad
 \text{in $K(G/\varpi^\Z\times D^\times/\varpi^\Z)$.}
 \]
\end{prop}

\begin{prf}
 We only give a sketch of a proof, since it seems to be well-known (\cf \cite[Proposition 4.4.13]{MR2074714}). 
 
 It is known that $M_\Dr/\varpi^\Z$ has an open covering $\mathfrak{U}=\{U_\lambda\}_{\lambda\in\Lambda}$
consisting of quasi-compact open subsets, indexed by the set $\Lambda$ of vertices of the Bruhat-Tits building
of $\mathrm{PGL}_n(F)$. This covering satisfies the following properties:
\begin{itemize}
 \item[(a)] For $g\in G/\varpi^\Z$, $U_\lambda\cdot g=U_{g^{-1}\lambda}$.
 \item[(b)] For each $\lambda\in\Lambda$, there exist only finitely many $\lambda'\in\Lambda$ satisfying $U_\lambda\cap U_{\lambda'}\neq\varnothing$. 
 \item[(c)] For each $\lambda\in\Lambda$, 
       $K_\lambda=\{g\in G/\varpi^\Z\mid U_\lambda\cdot g=U_\lambda\}$ is a compact open subgroup of $G/\varpi^\Z$.
\end{itemize}
For a finite subset $I\subset \Lambda$, put
\[
 U_I=\bigcap_{\lambda\in I}U_\lambda,\quad
 K_I=\{g\in G/\varpi^\Z\mid U_I\cdot g=U_I\}.
\]
For an integer $r\ge 0$, set $\Lambda_r=\{I\subset \Lambda\mid \#I=r+1, U_I\neq\varnothing\}$.
Then, by the three properties above, we have the following:
\begin{itemize}
 \item For each $r\ge 0$, the set of $G/\varpi^\Z$-orbits in $\Lambda_r$ is finite.
 \item We have $\Lambda_r=\varnothing$ for sufficiently large $r$.
 \item For every finite subset $I\subset \Lambda$, $K_I$ is a compact open subgroup of $G/\varpi^\Z$.
\end{itemize}
%
Take a system of representatives $I_{r,1},\ldots,I_{r,N_r}$ of $(G/\varpi^\Z)\backslash \Lambda_r$
and put $K_{r,i}=K_{I_{r,i}}$.
 
Let $\mathfrak{U}_m=\{U_{m,\lambda}\}_{\lambda\in\Lambda}$ be the covering
obtained as the inverse image of $\mathfrak{U}$. For a finite subset $I\subset\Lambda$,
put $U_{m,I}=\bigcap_{\lambda\in I}U_{m,\lambda}$. Then we have the \v{C}ech spectral sequence
 \[
 E_1^{-r,s}=\bigoplus_{I\in \Lambda_r}H^s_c(U_{m,I}\otimes_{\breve{F}}\overline{\breve{F}},\overline{\Q}_\ell)\Longrightarrow H^{-r+s}_{\Dr,m}.
 \]
 Put
 \[
 V^s_{m,r,i}=H^s_c(U_{m,I_{r,i}}\otimes_{\breve{F}}\overline{\breve{F}},\overline{\Q}_\ell).
 \]
 It is a finite-dimensional smooth representation of $K_{r,i}\times D^\times/\varpi^\Z$
 and vanishes for $s>2n$
 (\cf \cite[Propositon 5.5.1, Proposition 6.3.2]{MR1734903}).
 We can easily observe that $E_1^{-r,s}$ is isomorphic to $\bigoplus_{i=1}^{N_r}\cInd_{K_{r,i}}^{G/\varpi^\Z} V^s_{m,r,i}$ as a $G/\varpi^\Z\times D^\times/\varpi^\Z$-representation. 
 Therefore $E_1^{-r,s}$ is a finitely generated $G/\varpi^\Z$-representation,
 and vanishes for all but finitely many $(r,s)$. 
 Hence $H^i_{\Dr,m}$ is finitely generated as a $G$-module (\cf \cite[Remarque 3.12]{MR771671}).
 Moreover, in $K(G/\varpi^\Z\times D^\times/\varpi^\Z)$ we have
 \[
  \sum_i(-1)^i[H^i_{\Dr,m}]=\sum_{r,s}(-1)^{-r+s}E_1^{-r,s}=\sum_{r,s}\sum_{i=1}^{N_r}(-1)^{-r+s}\cInd_{K_{r,i}}^{G/\varpi^\Z} V^s_{m,r,i}.
 \]
This concludes the proof.
\end{prf}

\begin{defn}
 We denote by $\eta_m$ the image of $\sum_i (-1)^i[H^i_{\Dr,m}]$ under the rank map 
 \[
 \Rk\colon K(G/\varpi^\Z\times D^\times/\varpi^\Z)\longrightarrow \overline{\mathcal{H}}(G/\varpi^\Z\times D^\times/\varpi^\Z)
 \]
 (\cf \cite[1.2]{MR1779801}). For $h\in D^\times$, we define $\eta_{m,h}\in \overline{\mathcal{H}}(G/\varpi^\Z)$
 by
 \[
  \eta_{m,h}(g)=\int_{D^\times/\varpi^\Z}\eta_m(g,h'^{-1}hh')dh'.
 \]
\end{defn}

Using the expression of $\sum_i (-1)^i[H^i_{\Dr,m}]$ in Proposition \ref{prop:fin-gen},
we can give more explicit description of $\eta_m$ and $\eta_{m,h}$:

\begin{prop}\label{prop:eta-cInd}
 For $m\ge 0$ and $h\in D^\times$, define $\widetilde{\eta}_m\in \mathcal{H}(G/\varpi^\Z\times D^\times/\varpi^\Z)$
 and $\widetilde{\eta}_{m,h}\in \mathcal{H}(G/\varpi^\Z)$ by
 \[
 \widetilde{\eta}_m=\sum_{\nu=1}^N\frac{\varepsilon_\nu\theta_{\sigma_{m,\nu}^\vee}}{\vol(K_\nu\times D^\times/\varpi^\Z)},\qquad
 \widetilde{\eta}_{m,h}=\sum_{\nu=1}^N\frac{\varepsilon_\nu\theta_{\sigma_{m,\nu}^\vee}(-,h)}{\vol(K_\nu)},
 \]
 where $(-)^\vee$ denotes the contragredient, and 
 $\theta_{\sigma_{m,\nu}^\vee}$ (resp.\ $\theta_{\sigma_{m,\nu}^\vee}(-,h)$) is regarded as 
 a function on $G/\varpi^\Z\times D^\times/\varpi^\Z$ (resp.\ $G/\varpi^\Z$) by setting
 $\theta_{\sigma_{m,\nu}^\vee}(g,h)=0$ for $g\notin K_\nu$.

 Then, the image of $\widetilde{\eta}_m$ (resp.\ $\widetilde{\eta}_{m,h}$) in
 $\overline{\mathcal{H}}(G/\varpi^\Z\times D^\times/\varpi^\Z)$ (resp.\ $\overline{\mathcal{H}}(G/\varpi^\Z)$)
 coincides with $\eta_m$ (resp.\ $\eta_{m,h}$).
\end{prop}

\begin{prf}
 The assertion for $\widetilde{\eta}_{m,h}$ immediately follows from that for $\widetilde{\eta}_m$.
 Thus it suffices to prove the following:
 \begin{quote}
  Let $K$ be a compact open subgroup of $G/\varpi^\Z$.
  For every finite-dimensional smooth representation $\sigma$ of $K\times D^\times/\varpi^\Z$,
  the image of $\vol(K)^{-1}\theta_{\sigma^\vee}$ in $\overline{\mathcal{H}}(G/\varpi^\Z\times D^\times/\varpi^\Z)$
  coincides with $\Rk([\cInd_K^{G/\varpi^\Z} \sigma])$.
 \end{quote}
 Since the image of $[\sigma]$ under the rank map $\Rk\colon K(K\times D^\times/\varpi^\Z)\longrightarrow \overline{\mathcal{H}}(K\times D^\times/\varpi^\Z)$ is $\vol(K\times D^\times/\varpi^\Z)^{-1}\theta_{\sigma^\vee}$,
 this claim follows from the commutative diagram below (\cf \cite[proof of Theorem 1.6]{MR1779801}):
 \[
  \xymatrix{%
 K(K\times D^\times/\varpi^\Z)\ar[r]^-{\Rk}\ar[d]_-{\cInd_K^{G/\varpi^\Z}}& \overline{\mathcal{H}}(K\times D^\times/\varpi^\Z)\ar[d]^-{\text{\rm extention by $0$}}\\
 K(G/\varpi^\Z\times D^\times/\varpi^\Z)\ar[r]^-{\Rk}& \overline{\mathcal{H}}(G/\varpi^\Z\times D^\times/\varpi^\Z)\lefteqn{.}
 }
 \]
\end{prf}

In \cite{MR1302321}, Faltings investigated the function $\widetilde{\eta}_m$ above
by means of the Lefschetz trace formula. His results can be summarized in the following theorem.

\begin{thm}\label{thm:LTF}
 Let $g\in G^{\mathrm{reg}}$ and $h\in D^\times$.
 \begin{enumerate}
  \item If $g$ is elliptic, then we have
	\begin{align*}
	O^{G/\varpi^\Z}_g(\eta_{m,h})&=\#\Fix\bigl((g^{-1},h^{-1});M_{\Dr,m}/\varpi^\Z\bigr)\\
	 &=n\cdot \#\bigl\{a\in D^\times/\varpi^\Z(1+\Pi^m\mathcal{O}_D)\bigm| hah_g^{-1}=a\bigr\}.
	\end{align*}
  \item If $g$ is not elliptic, then we have $\int_{Z(g)\backslash G}\eta_{m,h}(x^{-1}gx)dx=0$, where
	$Z(g)$ denotes the centralizer of $g$.
 \end{enumerate}
\end{thm}

\begin{prf}
 Let us briefly recall the proof in \cite{MR1302321}.
 We use the notation in the proof of Proposition \ref{prop:fin-gen}.

 First consider the case where $g$ is elliptic. Then, we can find a finite subset $\Lambda_g\subset \Lambda$
 such that $g\Lambda_g=\Lambda_g$ and $gU_\lambda\cap U_\lambda=\varnothing$ for $\lambda\in \Lambda\setminus \Lambda_g$. Put $U_{m,g}=\bigcup_{\lambda\in \Lambda_g}U_{m,\lambda}$. Then $U_{m,g}$ is quasi-compact smooth and
 $(g^{-1},h^{-1})\colon U_{m,g}\longrightarrow U_{m,g}$ has no fixed point on the boundary of $U_{m,g}$.
 Therefore we can apply the Lefschetz trace formula for this endomorphism
 (for a general theory of the Lefschetz trace formula for rigid spaces, see \cite{adicLTF}).
 Noting that every fixed point of $(g^{-1},h^{-1})\colon M_{\Dr,m}/\varpi^\Z\longrightarrow M_{\Dr,m}/\varpi^\Z$
 lies in $U_{m,g}$, we obtain the following equality:
 \[
 \sum_i(-1)^i\Tr\bigl((g^{-1},h^{-1});H^i_c(U_{m,g}\otimes_{\breve{F}}\overline{\breve{F}},\overline{\Q}_\ell)\bigr)=\#\Fix\bigl((g^{-1},h^{-1});M_{\Dr,m}/\varpi^\Z\bigr).
 \]
 By the \v{C}ech spectral sequence, we can easily show that the left hand side is equal to
 $O^{G/\varpi^\Z}_g(\widetilde{\eta}_{m,h})$. The right hand side can be computed by using
 the period map, as in \cite[\S 2.6]{MR2383890}. The result is
 \[
  \#\Fix\bigl((g^{-1},h^{-1});M_{\Dr,m}/\varpi^\Z\bigr)=n\cdot \#\bigl\{a\in D^\times/\varpi^\Z(1+\Pi^m\mathcal{O}_D)\bigm| hah_g^{-1}=a\bigr\}.
 \]
 It is slightly different from \cite[Theorem 1]{MR1302321}, because our $M_{\Dr}/\varpi^\Z$ is
 the disjoint union of $n$ copies of $\Omega$ considered in \cite{MR1302321}.
 Rather, it is compatible with \cite[Theorem 2.6.8]{MR2383890}. This concludes the proof of i).

 To prove ii), apply the same argument to $M_{\Dr,m}/\Gamma$, where $\Gamma$ is a sufficiently small
 discrete torsion-free cocompact subgroup of $Z(g)$.
\end{prf}

\begin{cor}\label{cor:in-K0}
 For $\rho\in\Irr(D^\times/\varpi^\Z)$, $\Hom_{D^\times}(\rho,H^i_{\Dr})$
 is a finitely generated $G/\varpi^\Z$-representation by \cite[Lemma 5.2]{LT-LTF}.
 The image of $\sum_{i}(-1)^i[\Hom_{D^\times}(\rho,H^i_{\Dr})]$ under the map
 \[
 \Rk^\vee\colon K(G/\varpi^\Z)\yrightarrow{\Rk} \overline{\mathcal{H}}(G/\varpi^\Z)\yrightarrow{\vee} C^\infty(G^{\mathrm{ell}})
 \]
 coincides with $g\longmapsto n\theta_{\rho}(h_g^{-1})$.
 Recall that for $f\in \overline{\mathcal{H}}(G/\varpi^\Z)$
 the locally constant function $f^\vee$ on $G^{\mathrm{ell}}$ is given by
 $f^\vee(g)=\int_{G/\varpi^\Z}f(xg^{-1}x^{-1})dx$ (\cf \cite[p.~190]{MR1779801}).
\end{cor}

\begin{prf}
 Take a sufficiently large integer $m\ge 0$ so that $\rho\vert_{1+\Pi^m\mathcal{O}_D}$ is trivial.
 Then we have $\Hom_{D^\times}(\rho,H^i_{\Dr})=\Hom_{D^\times}(\rho,(H^i_{\Dr})^{1+\Pi^m\mathcal{O}_D})=\Hom_{D^\times}(\rho,H^i_{\Dr,m})$.

 It is easy to see that the following diagram is commutative:
 \[
  \xymatrix{%
 K(G/\varpi^\Z\times D^\times/\varpi^\Z)\ar[r]^-{\Rk}\ar[d]_-{\Hom_{D^\times}(\rho,-)}&
 \overline{\mathcal{H}}(G/\varpi^\Z\times D^\times/\varpi^\Z)\ar[d]^-{(*)}\\
 K(G/\varpi^\Z)\ar[r]^-{\Rk}&
 \overline{\mathcal{H}}(G/\varpi^\Z)\lefteqn{,}
 }
 \]
 where $(*)$ is given by
 \[
  f\longmapsto \Bigl(g\longmapsto \int_{D^\times/\varpi^\Z}f(g,h)\theta_\rho(h)dh\Bigr).
 \]
 Therefore, the image of $\sum_{i}(-1)^i[\Hom_{D^\times}(\rho,H^i_{\Dr,m})]$ under $\Rk^\vee$ can be calculated
 as follows:
 \begin{align*}
 g\longmapsto &\int_{G/\varpi^\Z}\int_{D^\times/\varpi^\Z}\eta_m(xg^{-1}x^{-1},h)\theta_\rho(h)dhdx\\
  &=\frac{1}{\vol(D^\times/\varpi^\Z)}\int_{G/\varpi^\Z}\int_{D^\times/\varpi^\Z}\int_{D^\times/\varpi^\Z}\eta_m(xg^{-1}x^{-1},h)\theta_\rho(h'hh'^{-1})dh'dhdx\\
  &=\frac{1}{\vol(D^\times/\varpi^\Z)}\int_{D^\times/\varpi^\Z}O^{G/\varpi^\Z}_{g^{-1}}(\eta_{m,h})\theta_{\rho}(h)dh\\
  &=\frac{n}{\vol(D^\times/\varpi^\Z)}\int_{D^\times/\varpi^\Z}\#\bigl\{a\in D^\times/\varpi^\Z(1+\Pi^m\mathcal{O}_D)\bigm| hah_g=a\bigr\}\theta_{\rho}(h)dh\\
  &=\frac{n}{\#(D^\times/\varpi^\Z(1+\Pi^m\mathcal{O}_D))}\\
  &\qquad\times\sum_{h\in D^\times/\varpi^\Z(1+\Pi^m\mathcal{O}_D)}
  \#\bigl\{a\in D^\times/\varpi^\Z(1+\Pi^m\mathcal{O}_D)\bigm| hah_g=a\bigr\}\theta_{\rho}(h)\\
  &=\frac{n}{\#(D^\times/\varpi^\Z(1+\Pi^m\mathcal{O}_D))}\sum_{a\in D^\times/\varpi^\Z(1+\Pi^m\mathcal{O}_D))}\theta_\rho(ah_g^{-1}a^{-1})\\
  &=n\theta_\rho(h_g^{-1}).
 \end{align*}
 This completes the proof.
\end{prf}

\section{Some harmonic analysis}\label{sec:harmonic-analysis}
In this section, we use Theorem \ref{thm:LTF} to investigate the virtual $D^\times$-representation 
$\sum_{i,j\ge 0}(-1)^{i+j}\Ext^j_{G/\varpi^\Z}(H^i_\Dr,\pi)$.
In the Lubin-Tate case, a similar study is carried out in \cite{LT-LTF}.
For $m\ge 0$, denote by $K'_m$ the image of $1+\Pi^m\mathcal{O}_D$ in $D^\times/\varpi^\Z$.

\begin{lem}\label{lem:Ext-invariant}
 For a smooth representation $V$ of $G/\varpi^\Z\times D^\times/\varpi^\Z$ and a smooth representation 
 $\pi$ of $G/\varpi^\Z$,
 we have $\Ext^j_{G/\varpi^\Z}(V,\pi)^{K'_m}\cong\Ext^j_{G/\varpi^\Z}(V^{K'_m},\pi)$.
 Here $\Ext_{G/\varpi^\Z}^j$ is taken in the category of smooth $G/\varpi^\Z$-representations.
\end{lem}

\begin{prf}
 First we prove that there exist a smooth $G/\varpi^\Z\times D^\times/\varpi^\Z$-representation $P$
 which is projective as a $G/\varpi^\Z$-representation and
 a $G/\varpi^\Z\times D^\times/\varpi^\Z$-equivariant surjection $P\longrightarrow V$.
 Take a $G/\varpi^\Z$-equivariant surjection $P'\longrightarrow V\longrightarrow 0$ from a projective
 $G/\varpi^\Z$-representation $P'$. Put $P=P'\otimes_\C C^\infty_c(D^\times/\varpi^\Z)$.
 Then $P$ is a smooth $G/\varpi^\Z\times D^\times/\varpi^\Z$-representation which is projective
 as a $G/\varpi^\Z$-representation, and a surjection $P\longrightarrow V$ is naturally induced.

 Therefore, we can take a $G/\varpi^\Z\times D^\times/\varpi^\Z$-equivariant resolution
 $P_\bullet\longrightarrow V\longrightarrow 0$ of $V$ such that $P_i$ is projective as a 
 smooth $G/\varpi^\Z$-representation. Since $P_i^{K'_m}$ is a direct summand of $P_i$
 as a $G/\varpi^\Z$-representation,
 $P_i^{K'_m}$ is also projective.
 Thus $P^{K'_m}_\bullet\longrightarrow V^{K'_m}\longrightarrow 0$
 gives a projective resolution of $V^{K'_m}$.
 Hence we have
 \begin{align*}
  \Ext^j_{G/\varpi^\Z}(V,\pi)^{K'_m}&=H^j\bigl(\Hom_{G/\varpi^\Z}(P_\bullet,\pi)\bigr)^{K'_m}
 \cong H^j\bigl(\Hom_{G/\varpi^\Z}(P^{K'_m}_\bullet,\pi)\bigr)\\
  &=\Ext^j_{G/\varpi^\Z}(V^{K'_m},\pi),
 \end{align*}
 as desired.
\end{prf}

\begin{cor}\label{cor:Ext-invariant-Dr}
 For every $m\ge 0$ and $\pi\in\Irr(G/\varpi^\Z)$,
 we have 
 \[
  \Ext^j_{G/\varpi^\Z}(H^i_{\Dr},\pi)^{K_m'}\cong \Ext^j_{G/\varpi^\Z}(H^i_{\Dr,m},\pi).
 \]
 It is finite-dimensional and vanishes if $j\ge n$.
 In particular, $\Ext^j_{G/\varpi^\Z}(H^i_{\Dr},\pi)^{\mathrm{sm}}$
 is an admissible representation of $D^\times/\varpi^\Z$ and vanishes if $j\ge n$,
 where $(-)^{\mathrm{sm}}$ denotes the set of $D^\times/\varpi^\Z$-smooth vectors.
\end{cor}

\begin{prf}
 Lemma \ref{lem:Ext-invariant} tells us that
 \[
  \Ext^j_{G/\varpi^\Z}(H^i_{\Dr},\pi)^{K_m'}
 =\Ext^j_{G/\varpi^\Z}\bigl((H^i_{\Dr})^{K_m'},\pi\bigr)
 =\Ext^j_{G/\varpi^\Z}(H^i_{\Dr,m},\pi).
 \]
 By Proposition \ref{prop:fin-gen} and \cite[Corollary II.3.3]{MR1471867},
 it is finite-dimensional and vanishes if $j\ge n$.
\end{prf}

\begin{rem}
 Later (Corollary \ref{cor:Bernstein}) we will prove that $\Ext^j_{G/\varpi^\Z}(H^i_{\Dr},\pi)$ is
 in fact a finite-dimensional smooth $D^\times$-representation.
\end{rem}

The character of $\sum_{i,j\ge 0}(-1)^{i+j}\Ext^j_{G/\varpi^\Z}(H^i_{\Dr,m},\pi)$ can be computed
by $\eta_{m,h}$ introduced in the previous section:

\begin{prop}\label{prop:Ext-char}
 For every $\pi\in\Irr(G/\varpi^\Z)$ and $h\in D^\times/\varpi^\Z$, we have
 \[
 \sum_{i,j\ge 0}(-1)^{i+j}\Tr\bigl(h;\Ext^j_{G/\varpi^\Z}(H^i_{\Dr,m},\pi)\bigr)
 =\Tr(\eta_{m,h};\pi)=\int_{G/\varpi^\Z}\eta_{m,h}(g)\theta_\pi(g)dg.
 \]
\end{prop}

\begin{prf}
 First note that $V\longmapsto \sum_{j\ge 0}(-1)^j\Tr(h,\Ext^j_{G/\varpi^\Z}(V,\pi))$ induces a
 homomorphism $K(G/\varpi^\Z)\longrightarrow \C$ of abelian groups. Therefore, by
 Proposition \ref{prop:fin-gen} and Proposition \ref{prop:eta-cInd}, we have only to show the following:
 \begin{quote}
  Let $K$ be a compact open subgroup of $G/\varpi^\Z$.
  For every finite-dimensional smooth representation $\sigma$ of $K\times D^\times/\varpi^\Z$ and $h\in D^\times$,
  we have
  \[
  \sum_{j\ge 0}(-1)^j\Tr\bigl(h;\Ext^j_{G/\varpi^\Z}(\cInd_K^{G/\varpi^\Z}\sigma,\pi)\bigr)=\Tr\Bigl(\frac{\theta_{\sigma^\vee}(-,h)}{\vol(K)};\pi\Bigr).
  \]
 \end{quote}
 Since $\cInd_K^{G/\varpi^\Z}\sigma$ is a projective $G/\varpi^\Z$-representation,
 $\Ext^j_{G/\varpi^\Z}(\cInd_K^{G/\varpi^\Z}\sigma,\pi)=0$ for $j\ge 1$. 
 Take an open normal subgroup $K_1\subset K$ such that $\sigma\vert_{K_1}$ is trivial.
 Then the left hand side can be
 computed as follows:
 \begin{align*}
  &\sum_{j\ge 0}(-1)^j\Tr\bigl(h;\Ext^j_{G/\varpi^\Z}(\cInd_K^{G/\varpi^\Z}\sigma,\pi)\bigr)
  =\Tr\bigl(h;\Hom_{G/\varpi^\Z}(\cInd_K^{G/\varpi^\Z}\sigma,\pi)\bigr)\\
  &\qquad =\Tr\bigl(h;\Hom_K(\sigma,\pi\vert_K)\bigr)=\Tr\bigl(h;\Hom_{K/K_1}(\sigma,\pi^{K_1})\bigr)\\
  &\qquad =\frac{1}{\#(K/K_1)}\sum_{g\in K/K_1}\Tr\bigl((g^{-1},h^{-1});\sigma\bigr)\Tr(g;\pi^{K_1})\\
  &\qquad =\frac{1}{\vol(K)}\Tr\bigl(\theta_{\sigma^\vee}(-,h);\pi\bigr).
 \end{align*}
 This completes the proof.
\end{prf}

\begin{lem}\label{lem:fix-orb-int}
 For every $g\in G^{\mathrm{ell}}$, $h\in D^\times$ and an integer $m\ge 0$, we have
  \[
  O^{G/\varpi^\Z}_g(\eta_{m,h})=nO^{D^\times/\varpi^\Z}_{h_g}\Bigl(\frac{\mathbf{1}_{hK'_m}}{\vol(K_m')}\Bigr),
 \]
 where $\mathbf{1}_{hK_m'}$ denotes the characteristic function of $hK_m'$.
\end{lem}

\begin{prf}
 By Theorem \ref{thm:LTF} i), we obtain
 \begin{align*}
  O^{G/\varpi^\Z}_g(\eta_{m,h})&=n\cdot \#\bigl\{a\in (D^\times/\varpi^\Z)/K'_m\bigm| hah_g^{-1}=a\bigr\}\\
 &=\frac{n}{\vol(K'_m)}\int_{D^\times/\varpi^\Z}\mathbf{1}_{hK'_m}(ah_ga^{-1})da
  =nO^{D^\times/\varpi^\Z}_{h_g}\Bigl(\frac{\mathbf{1}_{hK'_m}}{\vol(K'_m)}\Bigr).
 \end{align*}
\end{prf}

Next recall the definition of a transfer of a test function.

\begin{defn}
 For $\varphi\in \mathcal{H}(G)$ and $\varphi^D\in \mathcal{H}(D^\times)$,
 we say that $\varphi^D$ is a transfer of $\varphi$ if 
 \[
  \int_{D^\times/\varpi^\Z}\varphi^D(y^{-1}h_gy)dy=(-1)^{n-1}\int_{G/\varpi^\Z}\varphi(x^{-1}gx)dx
 \]
 for every $g\in G^{\mathrm{ell}}$.
\end{defn}

We know that if $\varphi\in \mathcal{H}(G)$ is supported on $G^{\mathrm{ell}}$, then
it has a transfer $\varphi^D\in \mathcal{H}(D^\times)$ (\cf \cite[Lemma 3.2]{LT-LTF}). 
The following lemma is obvious:

\begin{lem}\label{lem:phi-pi}
 Assume that $\varphi\in \mathcal{H}(G)$ is supported on $G^{\mathrm{ell}}$ and let
 $\varphi^D\in \mathcal{H}(D^\times)$ be its transfer.
 Put
 \[
  \varphi_{\varpi}(g)=\sum_{i\in\Z}\varphi(\varpi^ig),\qquad \varphi^D_{\varpi}(h)=\sum_{i\in\Z}\varphi^D(\varpi^ih).
 \]
 Then
 $\varphi_\varpi\in \mathcal{H}(G/\varpi^\Z)$, $\varphi^D_\varpi\in \mathcal{H}(D^\times/\varpi^\Z)$ and
 $O^{D^\times/\varpi^\Z}_{h_g}(\varphi^D_\varpi)=(-1)^{n-1}O^{G/\varpi^\Z}_g(\varphi_\varpi)$ 
 for every $g\in G^{\mathrm{ell}}$.
\end{lem}

\begin{thm}\label{thm:Dr-pi}
 Assume that $\varphi\in \mathcal{H}(G)$ is supported on $G^{\mathrm{ell}}$ and let
 $\varphi^D\in \mathcal{H}(D^\times)$  be its transfer.
 Then, for every $\pi\in \Irr(G/\varpi^\Z)$ we have
 \[
 \sum_{i,j\ge 0}(-1)^{i+j}\Tr\bigl(\varphi^D;\Ext^j_{G/\varpi^\Z}(H^i_\Dr,\pi)^{\mathrm{sm}}\bigr)
 =(-1)^{n-1}n\Tr(\varphi;\pi).
 \]
\end{thm}

\begin{prf}
 Let $\varphi_\varpi$ and $\varphi^D_\varpi$ be as in the previous lemma.
 Then clearly we have
 \begin{align*}
  \sum_{i,j\ge 0}(-1)^{i+j}\Tr\bigl(\varphi^D;\Ext^j_{G/\varpi^\Z}(H^i_\Dr,\pi)^{\mathrm{sm}}\bigr)
  &=\sum_{i,j\ge 0}(-1)^{i+j}\Tr\bigl(\varphi_\varpi^D;\Ext^j_{G/\varpi^\Z}(H^i_\Dr,\pi)^{\mathrm{sm}}\bigr),\\
  (-1)^{n-1}n\Tr(\varphi;\pi)&=(-1)^{n-1}n\Tr(\varphi_\varpi;\pi).
 \end{align*}
 Thus we may replace $\varphi$ and $\varphi^D$ by $\varphi_\varpi$ and $\varphi^D_\varpi$, respectively.

 Take $m\ge 0$ such that $\varphi^D_\varpi$ is $K'_m$-invariant, and
 write 
 \[
  \varphi^D_\varpi=\sum_{h\in J} a_h\frac{\mathbf{1}_{hK'_m}}{\vol(K'_m)},
 \]
 where $J$ is a finite subset of $D^\times/\varpi^\Z$ and $a_h\in\C$.
 By Corollary \ref{cor:Ext-invariant-Dr} and Proposition \ref{prop:Ext-char}, we have
 \begin{align*}
  \sum_{i,j\ge 0}(-1)^{i+j}\Tr\bigl(\varphi_\varpi^D;\Ext^j_{G/\varpi^\Z}(H^i_\Dr,\pi)^{\mathrm{sm}}\bigr)
  &=\sum_{i,j\ge 0,h\in J}(-1)^{i+j}a_h\Tr\bigl(h;\Ext^j_{G/\varpi^\Z}(H^i_{\Dr,m},\pi)\bigr)\\
  &=\sum_{h\in J}a_h\int_{G/\varpi^\Z}\eta_{m,h}(g)\theta_\pi(g)dg.
 \end{align*}
 By Theorem \ref{thm:LTF} ii) and Weyl's integral formula (\cf \cite[Theorem F]{MR874042}), we have
 \[
  \int_{G/\varpi^\Z}\eta_{m,h}(g)\theta_\pi(g)dg=\sum_{T}\frac{1}{\#W_T}\int_{T^{\mathrm{reg}}/\varpi^\Z}D(t)O^{G/\varpi^\Z}_t(\eta_{m,h})\theta_\pi(t)dt,
 \]
 where $T$ runs through conjugacy classes of elliptic maximal tori of $G$,
 $W_T$ denotes the rational Weyl group of $T$ and $D(t)$ denotes the Weyl denominator
 (\cf \cite[p.~185]{MR700135}). 
 The measure $dt$ on $T/\varpi^\Z$ is normalized so that the volume of $T/\varpi^\Z$ is one.
 Lemma \ref{lem:fix-orb-int} and Lemma \ref{lem:phi-pi} tell us that
 \begin{align*}
  \sum_{h\in J} a_hO^{G/\varpi^\Z}_t(\eta_{m,h})&=\sum_{h\in J} na_hO^{D^\times/\varpi^\Z}_{h_t}\Bigl(\frac{\mathbf{1}_{hK'_m}}{\vol(K_m')}\Bigr)=nO^{D^\times/\varpi^\Z}_{h_t}(\varphi_\varpi^D)\\
 &=(-1)^{n-1}nO^{G/\varpi^\Z}_t(\varphi_\varpi)
 \end{align*}
 for every $t\in T^{\mathrm{reg}}$. By Weyl's integral formula again, we have
 \begin{align*}
 &\sum_{i,j\ge 0}(-1)^{i+j}\Tr\bigl(\varphi^D_\varpi;\Ext^j_{G/\varpi^\Z}(H^i_\Dr,\pi)^{\mathrm{sm}}\bigr)\\
 &\qquad=(-1)^{n-1}n\sum_{T}\frac{1}{\#W_T}\int_{T^{\mathrm{reg}}/\varpi^\Z}D(t)O^{G/\varpi^\Z}_t(\varphi_\varpi)\theta_\pi(t)dt\\
  &\qquad=(-1)^{n-1}n\int_{G/\varpi^\Z}\varphi_\varpi(g)\theta_\pi(g)dg=(-1)^{n-1}n\Tr(\varphi_\varpi;\pi),
 \end{align*}
 as desired.
\end{prf}

\section{Faltings isomorphism}\label{sec:Faltings-isom}
Here we freely use the notation in \cite[Section 2]{LT-LTF}.
We need the following deep theorem due to Faltings (\cite{MR1936369}, see also \cite{MR2441311} for more
detailed exposition):

\begin{thm}\label{thm:Faltings}
 We have a $G\times D^\times$-equivariant isomorphism $H^i_{\mathrm{Dr}}\cong H^i_{\mathrm{LT}}$
 for every $i$.
\end{thm}

Note that the proof of Faltings' theorem does not require automorphic method.
It gives the following very important finiteness result on $H^i_{\mathrm{Dr}}$.

\begin{cor}\label{cor:Dr-adm}
 The $G$-representation $H^i_{\mathrm{Dr}}$ is admissible.
\end{cor}

\begin{prf}
 Put $K_m=\Ker(\GL_n(\mathcal{O})\to \GL_n(\mathcal{O}/(\varpi^m)))$. 
 Since 
 \[
  (H^i_{\mathrm{LT}})^{K_m}=H^i_c\bigl((M_m/\varpi^\Z)\otimes_{\breve{F}}\overline{\breve{F}},\overline{\Q}_\ell\bigr)
 \]
 is finite-dimensional, $H^i_{\mathrm{LT}}$ is an admissible representation of $G$.
 Thus, by Theorem \ref{thm:Faltings}, $H^i_{\mathrm{Dr}}$ is also an admissible representation of $G$.
\end{prf}

\begin{cor}\label{cor:Bernstein}
 For every $\pi\in \Irr(G/\varpi^\Z)$ and integers $i,j\ge 0$,
 $\Ext^j_{G/\varpi^\Z}(H^i_{\mathrm{Dr}},\pi)$ is a finite-dimensional smooth representation of $D^\times$.
 Moreover, $\Ext^j_{G/\varpi^\Z}(H^i_{\mathrm{Dr}},\pi)=0$ if $j\ge n$.
\end{cor}

\begin{prf}
 Let $\mathfrak{s}$ be the cuspidal support of $\pi$, and $H^i_{\mathrm{Dr},\mathfrak{s}}$ be
 the $\mathfrak{s}$-component of $H^i_{\mathrm{Dr}}$.
 Clearly we have $\Ext^j_{G/\varpi^\Z}(H^i_{\mathrm{Dr}},\pi)=\Ext^j_{G/\varpi^\Z}(H^i_{\mathrm{Dr},\mathfrak{s}},\pi)$.
 
 Let us observe that $H^i_{\mathrm{Dr},\mathfrak{s}}$ is a finitely generated $G$-representation.
 Since $H^i_{\mathrm{Dr},\mathfrak{s}}$ is an admissible $G$-representation by Corollary \ref{cor:Dr-adm},
 it is $\mathfrak{Z}(G)$-admissible, where $\mathfrak{Z}(G)$ denotes the Bernstein center of $G$
 (\cf \cite[\S 3.1]{MR771671}). Therefore, \cite[Corollaire 3.10]{MR771671} tells us that
 $H^i_{\mathrm{Dr},\mathfrak{s}}$ is finitely generated.
 In particular, there exists a compact open subgroup $K'\subset D^\times$ which acts
 on $H^i_{\mathrm{Dr},\mathfrak{s}}$ trivially. 

 Therefore, by \cite[Corollary II.3.3]{MR1471867}, $\Ext^j_{G/\varpi^\Z}(H^i_{\mathrm{Dr},\mathfrak{s}},\pi)$
 is finite-dimensional and vanishes if $j\ge n$.
 The natural action of $D^\times$ on $\Ext^j_{G/\varpi^\Z}(H^i_{\mathrm{Dr},\mathfrak{s}},\pi)$ is smooth,
 since the action of $K'$ is trivial. This completes the proof.
\end{prf}

\begin{defn}
 For $\pi\in\Irr(G/\varpi^\Z)$, put
 \[
  H_{\Dr}[\pi]=\sum_{i,j\ge 0}(-1)^{i+j}\Ext_{G/\varpi^\Z}^j(H^i_{\Dr},\pi)\quad \text{in $R(D^\times/\varpi^\Z)$.}
 \]
\end{defn}

We can consider the character $\theta_{H_{\Dr}[\pi]}$ of $H_{\Dr}[\pi]$.
Theorem \ref{thm:Dr-pi} can be written in the following way:

\begin{thm}\label{thm:char-rel-G}
 For every $\pi\in\Irr(G/\varpi^\Z)$ and $h\in (D^\times)^{\mathrm{reg}}$, we have
 \[
  \theta_{H_{\Dr}[\pi]}(h)=n\theta_\pi(g_h).
 \]
\end{thm}

\begin{prf}
 Theorem \ref{thm:Dr-pi} says that $\Tr(\varphi^D;H_{\Dr}[\pi])=(-1)^{n-1}n\Tr(\varphi;\pi)$.
 We can use exactly the same method as in the proof of \cite[Theorem 4.3]{LT-LTF}.
\end{prf}

The following is another consequence of Theorem \ref{thm:Faltings}:

\begin{cor}\label{cor:Dr-rho}
 For every $\rho\in\Irr(D^\times/\varpi^\Z)$, $H^i_{\Dr}[\rho]=\Hom_{D^\times}(H^i_{\Dr},\rho)^{\mathrm{sm}}$
 is a smooth $G$-representation of finite length. Moreover $H^i_{\Dr}[\rho]$ is isomorphic to
 $\Hom_{D^\times}(\rho,H^i_{\Dr})^\vee$.
\end{cor}

\begin{prf}
 By Proposition \ref{prop:fin-gen}, Corollary \ref{cor:Dr-adm} and \cite[Lemma 5.2]{LT-LTF}, 
 $H^i_{\Dr}[\rho]$ is a smooth $G$-representation of finite length.
 In the proof of \cite[Lemma 5.2]{LT-LTF}, a $G$-equivariant injection
 $H^i_{\Dr}[\rho]\hooklongrightarrow \Hom_{D^\times}(\rho,H^i_{\Dr})^\vee$ is constructed.
 It is easy to see that it is actually an isomorphism.
\end{prf}

\begin{defn}
 For $\rho\in\Irr(D^\times/\varpi^\Z)$, put
 \[
  H_{\Dr}[\rho]=\sum_{i}(-1)^iH^i_{\Dr}[\rho]\quad \text{in $R(G/\varpi^\Z)$.}
 \]
\end{defn}

We can consider the character $\theta_{H_{\Dr}[\rho]}$ of $H_{\Dr}[\rho]$.
Corollary \ref{cor:in-K0} can be written in the following way:

\begin{thm}\label{thm:char-rel-D}
 For every $\rho\in\Irr(D^\times/\varpi^\Z)$ and $h\in (D^\times)^{\mathrm{reg}}$, we have
 \[
  \theta_{H_{\Dr}[\rho]}(g_h)=n\theta_\rho(h).
 \]
\end{thm}

\begin{prf}
 We denote the natural homomorphism $R(G/\varpi^\Z)\longrightarrow K(G/\varpi^\Z)$ by $\mathrm{EP}$.
 By \cite[Lemma 3.7]{MR2308851}, the composite of
 \[
  R(G/\varpi^\Z)\yrightarrow{\mathrm{EP}} K(G/\varpi^\Z)\yrightarrow{\Rk}\overline{\mathcal{H}}(G/\varpi^\Z)\yrightarrow{\vee}C^\infty(G^{\mathrm{ell}})
 \]
 coincides with $\pi\longmapsto \theta_\pi\vert_{G^\mathrm{ell}}$  (it was originally proved in \cite[Theorem III.4.23]{MR1471867}). Therefore, by Corollary \ref{cor:in-K0} and Corollary \ref{cor:Dr-rho} we have
 \[
  \theta_{H_\Dr[\rho]}(g)=\bigl((\vee\circ \Rk\circ\mathrm{EP})(H_{\Dr[\rho]})\bigr)(g)=n\theta_\rho(h_g)
 \]
 for every $g\in G^{\mathrm{ell}}$. Hence $\theta_{H_{\Dr}[\rho]}(g_h)=n\theta_\rho(h)$ for every $h\in (D^\times)^\mathrm{reg}$, as desired.
\end{prf}

\begin{rem}
 The proof of \cite[Theorem III.4.23]{MR1471867} seems to use \cite[Theorem 0]{MR874042},
 whose proof relies on global technique. However, Theorem \ref{thm:Faltings} and
 \cite[Theorem 4.3]{LT-LTF} give an alternative proof of Theorem \ref{thm:char-rel-D},
 which does not involve any global argument.
\end{rem}

\section{Complements on representation theory}\label{sec:rep-theory}
For locally constant class functions $\varphi_1$, $\varphi_2$ on $G^{\mathrm{ell}}/\varpi^\Z$, put
\[
 \langle \varphi_1,\varphi_2\rangle_{\mathrm{ell}}
 =\sum_{T}\frac{1}{\#W_T}\int_{T/\varpi^\Z}D(t)\varphi_1(t)\overline{\varphi_2(t)}\,dt,
\]
where $T$ runs through conjugacy classes of elliptic maximal tori of $G$.
Other notations are also the same as in the proof of Theorem \ref{thm:Dr-pi}. 

For locally constant functions $\phi_1$, $\phi_2$ on $D^\times/\varpi^\Z$, put
\[
 \langle \phi_1,\phi_2\rangle
 =\int_{D^\times/\varpi^\Z}\phi_1(h)\overline{\phi_2(h)}\,dh,
\]
 where the measure $dh$ is normalized so that the volume of the compact group
 $D^\times/\varpi^\Z$ is one.
 
 These two pairings are compatible, in the sense of the following lemma:

\begin{lem}\label{lem:inner-prod}
 Let $\varphi_1$, $\varphi_2$ be locally constant class functions on $G^{\mathrm{ell}}/\varpi^\Z$, and 
 $\phi_1$, $\phi_2$ locally constant class functions on $D^\times/\varpi^\Z$.
 Assume that $\varphi_i(g_h)=\phi_i(h)$ for every $h\in (D^\times)^{\mathrm{reg}}$.
 Then, we have
 \[
 \langle \varphi_1,\varphi_2\rangle_{\mathrm{ell}}=\langle \phi_1,\phi_2\rangle.
 \]
\end{lem}

\begin{prf}
 Clear from Weyl's integral formula for $D^\times$.
\end{prf}

The following orthogonality relation of characters is very important for our work.

\begin{prop}\label{prop:orthogonality}
 For $\pi_1,\pi_2\in \Disc(G/\varpi^\Z)$, we have
 \[
  \langle \theta_{\pi_1},\theta_{\pi_2}\rangle_{\mathrm{ell}}=
 \begin{cases}1& \pi_1\cong \pi_2,\\ 0& \text{otherwise.}\end{cases}
 \]
\end{prop}

\begin{prf}
 Let $\omega_1$, $\omega_2$ be the central characters of $\pi_1$, $\pi_2$, respectively. 
 Then they are unitary, since $F^\times/\varpi^\Z$ is compact.
 If $\omega_1=\omega_2$, then the lemma follows immediately
 from \cite[Lemma 5.3]{MR700135}. Otherwise,
 \begin{align*}
  \int_{T/\varpi^\Z}D(t)\theta_{\pi_1}(t)\overline{\theta_{\pi_2}(t)}\,dt
  &=\int_{T/Z_G}\Bigl(\int_{Z_G/\varpi^\Z}D(tz)\theta_{\pi_1}(tz)\overline{\theta_{\pi_2}(tz)}\,dz\Bigr)dt\\
  &=\int_{T/Z_G}\Bigl(\int_{Z_G/\varpi^\Z}\omega_1(z)\overline{\omega_2(z)}\,dz\Bigr)D(t)\theta_{\pi_1}(t)\overline{\theta_{\pi_2}(t)}\,dt=0,
 \end{align*}
 as desired.
\end{prf}

\begin{rem}
 The proof of \cite[Theorem 5.3]{MR700135} (for example, \cite[\S A.3, \S A.4]{MR771672})
 seems to need a global argument, such as Howe's conjecture due to Clozel.
 However, at least if $n$ is prime, we can give a purely local proof of Proposition \ref{prop:orthogonality}
 as follows. Here we use freely the notation which will be introduced in the next section.

 Note that, if $n$ is prime, then any irreducible discrete series representation of $G$ is
 either a twisted Steinberg representation or supercuspidal (\cf \cite[Theorem 9.3]{MR584084}).
 First assume that both $\pi_1$ and $\pi_2$ are twisted Steinberg representations,
 and write $\pi_1=\mathbf{St}_{\chi_1}$ and $\pi_2=\mathbf{St}_{\chi_2}$,
 where $\chi_1$ and $\chi_2$ are characters of $F^\times/\varpi^\Z$.
 Then, by Lemma \ref{lem:inner-prod} and Lemma \ref{lem:Steinberg}, we have
 \[
  \langle \theta_{\mathbf{St}_{\chi_1}},\theta_{\mathbf{St}_{\chi_2}}\rangle_{\mathrm{ell}}
 =\langle \theta_{\chi_1\circ\Nrd},\theta_{\chi_2\circ\Nrd}\rangle
 =\int_{F^\times/\varpi^\Z}\chi_1(z)\overline{\chi_2(z)}dz
 =\begin{cases}1& \chi_1=\chi_2,\\ 0& \chi_1\neq \chi_2.\end{cases}
 \]
 Since $\chi_1\neq \chi_2$ implies $\mathbf{St}_{\chi_1}\ncong \mathbf{St}_{\chi_2}$
 (their characters are different), we have the orthogonality relation in this case.

 Next assume that $\pi_2$ is supercuspidal. Then, by \cite[\S A.3.e, \S A.3.g]{MR771672},
 we can find a matrix coefficient $\phi\in \mathcal{H}(G/\varpi^\Z)$ of $\pi_2$ satisfying the following:
 \begin{itemize}
  \item[(a)] $O_g^{G/\varpi^\Z}(\phi)=\overline{\theta_{\pi_2}(g)}$ for $g\in G^{\mathrm{ell}}$
	     and $\int_{Z(g)\backslash G}\phi(x^{-1}gx)dx=0$ for $g\in G^{\mathrm{reg}}\setminus G^{\mathrm{ell}}$.
  \item[(b)] $\Tr(\phi;\pi_2)=1$ and $\Tr(\phi;\pi_1)=0$ for $\pi_1\in\Disc(G/\varpi^\Z)$ with $\pi_1\ncong \pi_2$.
 \end{itemize}
 By (a) and Weyl's integral formula, we have
 \begin{align*}
 \Tr(\phi;\pi_1)&=\int_{G/\varpi^\Z} \phi(g)\theta_{\pi_1}(g)dg
  =\sum_T\frac{1}{\#W_T}\int_{T/\varpi^\Z}D(t)O_t^{G/\varpi^\Z}(\phi)\theta_{\pi_1}(t)dt\\
  &=\sum_T\frac{1}{\#W_T}\int_{T/\varpi^\Z}D(t)\overline{\theta_{\pi_2}(t)}\theta_{\pi_1}(t)dt
  =\langle\theta_{\pi_1},\theta_{\pi_2}\rangle_{\mathrm{ell}}.
 \end{align*}
 Therefore (b) gives the desired orthogonality relation.
\end{rem}

\section{Local Jacquet-Langlands correspondence}\label{sec:LJLC}
Let $R_I(G)$ be the submodule of $R(G)$ generated by the image of parabolically induced representations
(\cf \cite{MR874042}) and put $\overline{R}(G)=R(G)/R_I(G)$. It is known that 
$\overline{R}(G)$ is a free $\Z$-module with a basis $\{[\pi]\mid \pi\in\Disc(G)\}$
(\cf \cite[Lemme 2.1.4]{MR2308851}).
We regard $\Disc(G)$ as a subset of $\overline{R}(G)$.
Recall that the character of a parabolically induced representation vanishes on $G^{\mathrm{ell}}$.
Therefore, $\pi\mapsto \theta_\pi\vert_{G^{\mathrm{ell}}}$ induces a map
$\overline{R}(G)\longrightarrow C^\infty(G^{\mathrm{ell}})$.

Moreover, we denote by $\overline{R}(G/\varpi^\Z)$ the image of $R(G/\varpi^\Z)$ in $\overline{R}(G)$.
It is easy to see that $\{[\pi]\mid \pi\in\Disc(G/\varpi^\Z)\}$ gives a basis of $\overline{R}(G/\varpi^\Z)$.
Set $\overline{R}(G/\varpi^\Z)_\Q=\overline{R}(G/\varpi^\Z)\otimes_\Z\Q$
and $\overline{R}(G)_\Q=\overline{R}(G)\otimes_\Z\Q$.

\begin{lem}
 The homomorphism $\widetilde{\LJ}\colon R(G/\varpi^\Z)\longrightarrow R(D^\times/\varpi^\Z)_\Q$ given by 
 \[
 \pi\longmapsto \frac{(-1)^{n-1}}{n}H_\Dr[\pi]\quad \text{for $\pi\in\Irr(G/\varpi^\Z)$} 
 \]
 factors through $\overline{R}(G/\varpi^\Z)$.
\end{lem}

\begin{prf}
 Let $P$ be a proper parabolic subgroup of $G$ with Levi factor $M$
 and $\sigma$ an irreducible smooth representation of $M$ on which $\varpi\in Z_M$ acts trivially.
 By Theorem \ref{thm:char-rel-G}, the character of $\widetilde{\LJ}(\Ind_P^G\sigma)$
 on $(D^\times)^\mathrm{reg}$ is given by $h\longmapsto (-1)^{n-1}\theta_{\Ind_P^G\sigma}(g_h)=0$.
 Since $(D^\times)^{\mathrm{reg}}$ is dense in $D^\times$, it vanishes for every $h\in D^\times$.
 Thus $\widetilde{\LJ}(\Ind_P^G\sigma)=0$ by linear independence of characters.
 Since the kernel of $R(G/\varpi^\Z)\longrightarrow \overline{R}(G/\varpi^\Z)$
 is generated by such representations as $\Ind_P^G\sigma$, we conclude the proof.
\end{prf}

The following is the main construction in this paper.

\begin{defn}
 We define two homomorphisms 
 \[
  \JL\colon R(D^\times/\varpi^\Z)_\Q\longrightarrow \overline{R}(G/\varpi^\Z)_\Q,\qquad
  \LJ\colon \overline{R}(G/\varpi^\Z)_\Q\longrightarrow R(D^\times/\varpi^\Z)_\Q
 \] 
  by 
 \[
  \JL(\rho)=\frac{(-1)^{n-1}}{n}H_{\Dr}[\rho],\qquad
  \LJ(\pi)=\frac{(-1)^{n-1}}{n}H_{\Dr}[\pi]
 \]
 for $\rho\in\Irr(D^\times/\varpi^\Z)$ and $\pi\in\Irr(G/\varpi^\Z)$.
 The latter map is well-defined by the previous lemma.
\end{defn}

\begin{prop}\label{prop:JL-LJ-properties}
 \begin{enumerate}
  \item We have the character relations
	\[
	 \theta_\rho(h)=(-1)^{n-1}\theta_{\JL(\rho)}(g_h),\quad 
	\theta_\pi(g_h)=(-1)^{n-1}\theta_{\LJ(\pi)}(h)
	\]
	for every $\rho\in\Irr(D^\times/\varpi^\Z)$, $\pi\in\Irr(G/\varpi^\Z)$ and $h\in (D^\times)^\mathrm{reg}$.
  \item For every $\rho,\rho'\in \Irr(D^\times/\varpi^\Z)$ and $\pi,\pi'\in \Irr(G/\varpi^\Z)$,
	we have 
	\begin{gather*}
	 \langle\theta_{\JL(\rho)},\theta_{\JL(\rho')}\rangle_{\mathrm{ell}}
	=\langle\theta_{\rho},\theta_{\rho'}\rangle,\quad
	\langle\theta_{\LJ(\pi)},\theta_{\LJ(\pi')}\rangle
	=\langle\theta_{\pi},\theta_{\pi'}\rangle_{\mathrm{ell}},\\
	\langle\theta_{\JL(\rho)},\theta_\pi\rangle_{\mathrm{ell}}
	=\langle\theta_{\rho},\theta_{\LJ(\pi)}\rangle.
	\end{gather*}
  \item Two maps $\JL$ and $\LJ$ are inverse to each other. 
  \item The map $\JL$ is compatible with character twists.
	Namely, for a character $\chi$ of $F^\times$ which is trivial on $\varpi^{n\Z}\subset F^\times$,
	we have $\JL(\rho\otimes(\chi\circ\Nrd))=\JL(\rho)\otimes (\chi\circ\det)$.
	The same holds for $\LJ$.
  \item The map $\JL$ preserves central characters.
	Namely, for $\rho\in \Irr(D^\times/\varpi^\Z)$, 
	write $\JL(\rho)=\sum_{\pi\in\Disc(G/\varpi^\Z)}a_\pi[\pi]$.
	Then, every $\pi$ with $a_\pi\neq 0$ has the same central character as $\rho$.
	The same holds for $\LJ$.
 \end{enumerate}
\end{prop}

\begin{prf}
 i) is clear from Theorem \ref{thm:char-rel-G} and Theorem \ref{thm:char-rel-D}.
 ii) follows from i) and Lemma \ref{lem:inner-prod}.

 Prove iii). For $\pi\in\Disc(G/\varpi^\Z)$, write $\JL(\LJ(\pi))=\sum_{\pi'\in\Disc(G/\varpi^\Z)}a_{\pi'}[\pi']$.
 Then, by ii) and Proposition \ref{prop:orthogonality} we have
 \[
  a_{\pi'}=\langle\theta_{\JL(\LJ(\pi))},\theta_{\pi'}\rangle_{\mathrm{ell}}
 =\langle \theta_{\LJ(\pi)},\theta_{\LJ(\pi')}\rangle
 =\langle \theta_{\pi},\theta_{\pi'}\rangle_{\mathrm{ell}}.
 \]
 Therefore, $a_{\pi'}=1$ if $\pi'=\pi$, and $a_{\pi'}=0$ otherwise.
 In other words, $\JL(\LJ(\pi))=\pi$. Thus we have $\JL\circ\LJ=\id$.
 Similarly we can prove that $\LJ\circ\JL=\id$.

 For iv), it suffices to show that $\LJ$ is compatible with character twists.
 Let $\pi\in\Irr(G/\varpi^\Z)$ and $\chi$ be a character of $F^\times$ which is trivial on
 $\varpi^{n\Z}\subset F^\times$. Then, for every $h\in (D^\times)^{\mathrm{reg}}$, we have
 \begin{align*}
 \theta_{\LJ(\pi\otimes (\chi\circ \det))}(h)&=(-1)^{n-1}\theta_{\pi\otimes(\chi\circ \det)}(g_h)
 =(-1)^{n-1}\chi(\det g_h)\theta_\pi(g_h)\\
  &=\chi(\Nrd h)\theta_{\LJ(\pi)}(h)=\theta_{LJ(\pi)\otimes (\chi\circ\Nrd)}(h).
 \end{align*}
 Since $(D^\times)^{\mathrm{reg}}$ is dense in $D^\times$, we have 
 $\theta_{\LJ(\pi\otimes (\chi\circ \det))}=\theta_{LJ(\pi)\otimes (\chi\circ\Nrd)}$.
 By linear independence of characters, 
 we conclude that $\LJ(\pi\otimes (\chi\circ \det))=LJ(\pi)\otimes (\chi\circ\Nrd)$.

 Finally we prove v). Write $\JL(\rho)=\sum_{\pi\in\Disc(G/\varpi^\Z)}a_\pi[\pi]$.
 By Proposition \ref{prop:orthogonality}, $a_\pi=\langle \theta_{\JL(\rho)},\theta_{\pi}\rangle_{\mathrm{ell}}$.
 Assume that the central character of $\pi\in \Disc(G/\varpi^\Z)$ is different from
 that of $\rho$. Then, we have
 \[
  a_\pi=\langle \theta_{\JL(\rho)},\theta_{\pi}\rangle_{\mathrm{ell}}
 =\sum_T\frac{1}{\#W_T}\int_{T/\varpi^\Z}D(t)\theta_{\JL(\rho)}(t)\overline{\theta_{\pi}(t)}dt=0
 \]
 in the same way as in the proof of Proposition \ref{prop:orthogonality}.
 Therefore $\JL$ preserves central characters.
 By this result and ii), we have $\langle \theta_\rho,\theta_{\LJ(\pi)}\rangle=0$ unless $\rho$ and $\pi$ have
 the same central character. This means that the coefficient of $\rho\in \Irr(D^\times/\varpi^\Z)$
 in $\LJ(\pi)$ is zero unless the central character of $\rho$ is the same as that of $\pi$.
 Namely, $\LJ$ preserves central characters.
\end{prf}

By twisting, we can extend $\JL$ and $\LJ$ to maps
between $R(D^\times)$ and $\overline{R}(G)$.

\begin{prop}
 There exists a unique extension of $\JL$ to a homomorphism from
 $R(D^\times)_\Q$ to $\overline{R}(G)_\Q$ which is compatible with character twists. 
 Similarly, we have a unique extension of $\LJ$ to a homomorphism $\overline{R}(G)_\Q\longrightarrow R(D^\times)_\Q$
 which is compatible with character twists.
 We denote them by $\JL$ and $\LJ$ again.
 These are inverse to each other, satisfy the same character relations as in 
 Proposition \ref{prop:JL-LJ-properties} i), and preserve central characters.
\end{prop}

\begin{prf}
 For $\rho\in \Irr(D^\times)$, let $\omega_\rho$ be its central character.
 Take $c\in \C^\times$ such that $c^n=\omega_{\rho}(\varpi)$, and consider the character
 $\chi_c\colon z\longmapsto c^{v_F(z)}$ of $F^\times$.
 Then $\rho\otimes (\chi_c^{-1}\circ\Nrd)\in \Irr(D^\times/\varpi^\Z)$.
 Extend $\JL$ to $R(D^\times)_\Q\longrightarrow \overline{R}(G)_\Q$ by 
 \[
  \JL(\rho)=\JL\bigl(\rho\otimes (\chi_c^{-1}\circ \Nrd)\bigr)\otimes (\chi_c\circ \det).
 \]
 By Proposition \ref{prop:JL-LJ-properties} iv), it is independent of the choice of $c$.
 Moreover, we can easily observe that it is the unique extension of the original $\JL$
 which is compatible with character twists.
 Similarly, we can uniquely extend $\LJ$ to a map $\overline{R}(G)_\Q\longrightarrow R(D^\times)_\Q$
 compatible with character twists.

 By using Proposition \ref{prop:JL-LJ-properties} v), we can easily check that the extended $\JL$ and $\LJ$ are
 inverse to each other.
 The remaining parts are also immediate consequences of Proposition \ref{prop:JL-LJ-properties} i), v).
\end{prf}

Next we will observe the uniqueness of the maps $\JL$, $\LJ$ satisfying the character relations.

\begin{prop}
 \begin{enumerate}
  \item Let $\JL'\colon R(D^\times)_\Q\longrightarrow \overline{R}(G)_\Q$ be a homomorphism
	satisfying the character relation
	$\theta_\rho(h)=(-1)^{n-1}\theta_{\JL'(\rho)}(g_h)$ for every $h\in (D^\times)^{\mathrm{reg}}$.
	Then we have $\JL'=\JL$.
  \item Let $\LJ'\colon \overline{R}(G)_\Q\longrightarrow R(D^\times)_\Q$ be a homomorphism
	satisfying the character relation
	$\theta_\pi(g_h)=(-1)^{n-1}\theta_{\LJ'(\pi)}(h)$ for every $h\in (D^\times)^{\mathrm{reg}}$.
	Then we have $\LJ'=\LJ$.
 \end{enumerate}
\end{prop}

\begin{prf}
 To prove i), it suffices to show that $\LJ\circ \JL'=\id$.
 By the character relation, we have $\theta_{\LJ(\JL'(\rho))}(h)=\theta_\rho(h)$
 for every $\rho\in\Irr(D^\times)$ and $h\in (D^\times)^\mathrm{reg}$.
 Thus we can conclude that $\LJ(\JL'(\rho))=\rho$ by linear independence of characters for $D^\times$.
 For ii), prove $\LJ'\circ \JL=\id$ by a similar argument.
\end{prf}

So far, we have obtained the following theorem:

\begin{thm}\label{thm:main}
 We can construct the following two homomorphisms geometrically:
 \[
  \JL\colon R(D^\times)_\Q\longrightarrow \overline{R}(G)_\Q,
 \quad
 \LJ\colon \overline{R}(G)_\Q\longrightarrow R(D^\times)_\Q.
 \]
 These two maps are inverse to each other, and satisfy
 the character relations
 \[
  \theta_\rho(h)=(-1)^{n-1}\theta_{\JL(\rho)}(g_h),\quad 
 \theta_\pi(g_h)=(-1)^{n-1}\theta_{\LJ(\pi)}(h)
 \] 
 for every $h\in (D^\times)^\mathrm{reg}$.
 They are characterized by these character relations.
 Moreover, $\JL$ and $\LJ$ are compatible with character twists, and preserve central characters.
\end{thm}

Let $B\subset G$ be the Borel subgroup consisting of upper triangular matrices.
Recall that the Steinberg representation $\mathbf{St}$ is the unique irreducible quotient of 
the unnormalized induction $\Ind_B^G\mathbf{1}$ from the trivial character $\mathbf{1}$ on $B$.
For a character $\chi$ of $F^\times$, put
$\mathbf{St}_\chi=\mathbf{St}\otimes(\chi\circ \det)$. A representation of the form $\mathbf{St}_\chi$
is called a twisted Steinberg representation. It is an irreducible discrete series representation of $G$.
The following lemma is very well-known:

\begin{lem}\label{lem:Steinberg}
 We have
 $\theta_{\chi\circ\Nrd}(h)=(-1)^{n-1}\theta_{\mathbf{St}_{\chi}}(g_h)$
 for a character $\chi$ of $F^\times$ and $h\in (D^\times)^{\mathrm{reg}}$.
\end{lem}

\begin{prf}
 In $\overline{R}(G)$, we have $[\mathbf{St}_{\chi}]=(-1)^{n-1}[\chi\circ\det]$
 (\cf \cite[Remarque 2.1.14]{MR2308851}).
 As the character of a parabolically induced representation vanishes on $G^{\mathrm{ell}}$,
 we have
 \begin{align*}
  \theta_{\mathbf{St}_\chi}(g_h)&=(-1)^{n-1}\theta_{\chi\circ\det}(g_h)=(-1)^{n-1}\chi(\det g_h)=(-1)^{n-1}\chi(\Nrd h)\\
  &=(-1)^{n-1}\theta_{\chi\circ\Nrd}(h),
 \end{align*}
 as desired.
\end{prf}

\begin{cor}\label{cor:Steinberg}
 For a character $\chi$ of $F^\times$, we have $\JL(\chi\circ\Nrd)=\mathbf{St}_\chi$ and
 $\LJ(\mathbf{St}_\chi)=\chi\circ\Nrd$.
\end{cor}

\begin{prf}
 By Lemma \ref{lem:Steinberg}, we have $\theta_{\LJ(\mathbf{St}_\chi)}=\theta_{\chi\circ\Nrd}$.
 Linear independence of characters tells us that $\LJ(\mathbf{St}_\chi)=\chi\circ\Nrd$.
\end{prf}

The following is a consequence of the non-cuspidality result in \cite{non-cusp}:

\begin{prop}\label{prop:non-cusp}
 For an irreducible supercuspidal representation $\pi$ of $G$,
 write $\LJ(\pi)=\sum_{\rho\in\Irr(D^\times)}a_\rho[\rho]$. Then we have $a_\rho\ge 0$ for every $\rho$.
\end{prop}

\begin{prf}
 We may assume that $\pi\in \Irr(G/\varpi^\Z)$.
 Since $\pi$ is injective in the category of smooth $G/\varpi^\Z$-representations,
 $\Ext^j_{G/\varpi^\Z}(H_\Dr^i,\pi)=0$ unless $j=0$.
 By Theorem \ref{thm:Faltings} and \cite[Theorem 3.7]{non-cusp}, we have
 $\Hom_{G/\varpi^\Z}(H_\Dr^i,\pi)=0$ unless $i=n-1$.
 Therefore we have $\LJ(\pi)=n^{-1}[\Hom_{G/\varpi^\Z}(H_\Dr^{n-1},\pi)]$.
 This concludes the proof.
\end{prf}

Now we can prove the local Jacquet-Langlands correspondence for prime $n$.

\begin{thm}\label{thm:LJLC-prime}
 Assume that $n$ is a prime number.
 Then $\JL$ induces a bijection 
 \[
  \JL\colon \Irr(D^\times)\yrightarrow{\cong} \Disc(G)
 \]
 satisfying the character relation
 $\theta_\rho(h)=(-1)^{n-1}\theta_{\JL(\rho)}(g_h)$
 for every $h\in (D^\times)^{\mathrm{reg}}$.
\end{thm}

\begin{prf}
 For simplicity, we denote by $\Cusp(G/\varpi^\Z)$ the subset of $\Disc(G/\varpi^\Z)$
 consisting of supercuspidal representations.
 By Theorem \ref{thm:main}, it suffices to show the following:
 \begin{itemize}
  \item[(a)] For $\rho\in\Irr(D^\times/\varpi^\Z)$, $\JL(\rho)\in\Disc(G/\varpi^\Z)$.
  \item[(b)] For $\pi\in\Disc(G/\varpi^\Z)$, $\LJ(\pi)\in\Irr(D^\times/\varpi^\Z)$.
 \end{itemize}
 First we shall prove (a). If $\rho$ is a character, then it follows from Corollary \ref{cor:Steinberg}.
 Assume that $\rho$ is not a character, and write $\JL(\rho)=\sum_{\pi\in\Disc(G/\varpi^\Z)}a_\pi[\pi]$.
 Then, $a_\pi=0$ unless $\pi$ is supercuspidal. Indeed, if $\pi\notin \Cusp(G/\varpi^\Z)$,
 then $\pi$ is a twisted Steinberg representation $\mathbf{St}_\chi$, for $n$ is prime 
 (\cf \cite[Theorem 9.3]{MR584084}).
 By Proposition \ref{prop:orthogonality}, Proposition \ref{prop:JL-LJ-properties} ii) and
 Corollary \ref{cor:Steinberg}, we have
 \[
  a_\pi=a_{\mathbf{St}_\chi}=\langle \theta_{\JL(\rho)},\theta_{\mathbf{St}_\chi}\rangle_{\mathrm{ell}}
 =\langle \theta_{\rho},\theta_{\LJ(\mathbf{St}_\chi)}\rangle
 =\langle \theta_{\rho},\theta_{\chi\circ\Nrd}\rangle=0.
 \]
 Since $\JL(\rho)\neq 0$, there exists at least one $\pi\in \Cusp(G/\varpi^\Z)$
 satisfying $a_\pi\neq 0$. 
 Let us observe that such $\pi$ is unique. Assume that there exist $\pi, \pi'\in \Cusp(G/\varpi^\Z)$
 such that $a_\pi$ and $a_{\pi'}$ are non-zero.
 Then, we have
 $\langle \theta_{\LJ(\pi)},\theta_\rho\rangle=\langle \theta_{\pi},\theta_{\JL(\rho)}\rangle_{\mathrm{ell}}=a_\pi\neq 0$, and similarly $\langle \theta_{\LJ(\pi')},\theta_\rho\rangle\neq 0$.
 In other words, if we write 
 \[
  \LJ(\pi)=\sum_{\varrho\in\Irr(D^\times/\varpi^\Z)}b_\varrho[\varrho],\quad
 \LJ(\pi')=\sum_{\varrho\in\Irr(D^\times/\varpi^\Z)}b'_\varrho[\varrho],
 \]
 then $b_\rho=a_\pi$ and $b'_\rho=a_{\pi'}$ are non-zero. On the other hand, Proposition \ref{prop:non-cusp}
 tells us that $b_\varrho\ge 0$ and $b'_\varrho\ge 0$ for every $\varrho$.
 Thus, by Proposition \ref{prop:JL-LJ-properties} ii), we conclude that 
 \[
 \langle \theta_\pi,\theta_{\pi'}\rangle_{\mathrm{ell}}=\langle \theta_{\LJ(\pi)},\theta_{\LJ(\pi')}\rangle
 =\sum_{\varrho}b_\varrho b'_\varrho>0,
 \]
 which is equivalent to $\pi\cong \pi'$ by Proposition \ref{prop:orthogonality}.
 Now we have $\JL(\rho)=a_\pi[\pi]$ for some $\pi\in\Cusp(G/\varpi^\Z)$.
 Moreover, the argument above tells us that $a_\pi\ge 0$.
 By Proposition \ref{prop:orthogonality} and Proposition \ref{prop:JL-LJ-properties} ii), we have
 \[
 1=\langle\theta_{\rho},\theta_{\rho}\rangle=\langle\theta_{\JL(\rho)},\theta_{\JL(\rho)}\rangle_{\mathrm{ell}}=a_\pi^2\langle\theta_\pi,\theta_\pi\rangle_{\mathrm{ell}}=a^2_\pi.
 \]
 Hence we conclude that $a_\pi=1$ and $\JL(\rho)=\pi\in \Disc(G/\varpi^\Z)$.

 Next prove (b). Write $\LJ(\pi)=\sum_{\rho\in\Irr(D^\times/\varpi^\Z)}b_\rho[\rho]$.
 Then, since $\JL$ and $\LJ$ are inverse to each other,
 we have $\pi=\sum_{\rho\in\Irr(D^\times/\varpi^\Z)}b_\rho\JL(\rho)$, and thus
 \[
  1=\sum_{\rho\in\Irr(D^\times/\varpi^\Z)}b_\rho\langle \theta_{\JL(\rho)},\theta_\pi\rangle_{\mathrm{ell}}.
 \]
 By (a) and Proposition \ref{prop:orthogonality},
 there exists $\rho\in \Irr(D^\times/\varpi^\Z)$ such that $\pi=\JL(\rho)$.
 Then $\LJ(\pi)=\rho\in \Irr(D^\times/\varpi^\Z)$, as desired.
\end{prf}

\def\cftil#1{\ifmmode\setbox7\hbox{$\accent"5E#1$}\else
  \setbox7\hbox{\accent"5E#1}\penalty 10000\relax\fi\raise 1\ht7
  \hbox{\lower1.15ex\hbox to 1\wd7{\hss\accent"7E\hss}}\penalty 10000
  \hskip-1\wd7\penalty 10000\box7} \def\cprime{$'$} \def\cprime{$'$}
\providecommand{\bysame}{\leavevmode\hbox to3em{\hrulefill}\thinspace}
\providecommand{\MR}{\relax\ifhmode\unskip\space\fi MR }
\providecommand{\MRhref}[2]{%
  \href{http://www.ams.org/mathscinet-getitem?mr=#1}{#2}
}
\providecommand{\href}[2]{#2}

\end{document}